\documentclass[a4paper,10pt]{article}

\usepackage{t1enc}
\usepackage[latin1]{inputenc}
\usepackage{graphics}
\usepackage{graphicx}
\usepackage{amsfonts}
\usepackage{eufrak}
\usepackage{xypic}
\usepackage{amssymb}
\usepackage{amsmath}
\usepackage{amsthm}
\usepackage{mathrsfs}

\theoremstyle{plain}
\newtheorem{thm}{Theorem}[section]
\newtheorem{lemma}[thm]{Lemma}
\newtheorem{prop}[thm]{Proposition}
\newtheorem{cor}[thm]{Corollary}

\theoremstyle{definition}
\newtheorem{df}[thm]{Definition}
\newtheorem{rem}[thm]{Remark}
\newtheorem{ex}[thm]{Example}

\title{\bf On the growth of Hermitian groups}
\author{Rui Palma\\  {\footnotesize\emph{E-mail: ruip@math.uio.no}}}
\date{}

\begin{document}

\maketitle

\begin{abstract}
A locally compact group $G$ is said to be Hermitian if every selfadjoint element of $L^1(G)$ has real spectrum. Using Halmos' notion of capacity in Banach algebras and a result of Jenkins, Fountain, Ramsay and Williamson we will put a bound on the growth of Hermitian groups. In other words, we will show that if $G$ has a subset that grows faster than a certain constant, then $G$ cannot be Hermitian. Our result allows us to give new examples of non-Hermitian groups which could not tackled by the existing theory. The examples include certain infinite free Burnside groups, automorphism groups of trees, and $p$-adic general and special linear groups.
\end{abstract}

{\renewcommand{\thefootnote}{}
\footnotetext{\emph{Date:} \today}}

\section*{Introduction}

A locally compact group $G$ is said to be \emph{Hermitian} when $L^1(G)$ is a Hermitian Banach $^*$-algebra, i.e. when every selfadjoint element of $L^1(G)$ has real spectrum. There are many classes of locally compact groups which are known to be Hermitian, and these include abelian groups, compact groups, nilpotent groups, $FC^{-}$-groups and also a wide class of Lie groups (the reader is referred to \cite[12.6.22]{palmer} for an account).

The class of Hermitian groups has been more successfully studied in the case of connected locally compact groups. In this setting it is known, by a result of Jenkins \cite[Theorem 4.5]{jenkins nonsym}, that a connected, reductive Lie group is Hermitian if and only if its semisimple quotient is compact. Moreover, it was shown by Palmer \cite[Theorem 12.5.18 (e)]{palmer} that every almost connected Hermitian locally compact group is necessarily amenable. Both these results automatically provide us with many examples of non-Hermitian groups, as for instance, any non-amenable connected Lie group.

In the case of discrete groups the most important result is due to Jenkins \cite[Theorem 5.1]{jenkins}, and states that a discrete Hermitian group cannot contain a free sub-semigroup in two generators. This result also automatically provides us with many examples of non-Hermitian discrete groups, such as all non-abelian free groups and all solvable groups of exponential growth, for example.

The question of Hermitianess in the case of non-discrete totally disconnected groups has been, on the other hand, largely unaddressed in the literature, with no general results known and with several concrete examples still to be decided whether they are Hermitian or not. According to Palmer \cite[Table 5, pages 1488-1490]{palmer} such groups, for which the question remains unanswered, include the automorphism groups of trees (in their totally disconnected topology) and the $p$-adic "ax+b" group $\mathbb{Q}_p \rtimes \mathbb{Q}_p^*$. In fact, we only know of one example of a non-discrete totally disconnected group which has been proven to be non-Hermitian, and that is $PGL_2(\mathbb{Q}_p)$, as shown by Jenkins in \cite[page 300]{jenkins nonsym}.

The goal of this article is to establish a connection between growth rates in locally compact groups and Hermitianess. Our main result says, essentially, that if a locally compact group has a subset whose growth rate is larger than a certain fixed number, than the group is not Hermitian. In other words, subsets of Hermitian locally compact groups cannot grow very fast.

Our result was inspired and based on the, seemingly independent, works of Jenkins \cite{jenkins nonsym} and Fountain, Ramsay and Williamson \cite{fountain}. Using Halmos' notion of \emph{capacity} in Banach algebras \cite{halmos}, they gave a sufficient condition for a Banach $^*$-algebra to be non-Hermitian. This allowed Fountain, Ramsay and Williamson to give an alternative proof that the free group on two generators is not Hermitian based on the very fast growth of this group. Our goal is to extend this idea to more general locally compact groups.

As a consequence of our result we can give new examples of non-Hermitian groups which could not be tackled by the existing theory. For instance, we can show that certain free Burnside groups are not Hermitian, being the first examples of discrete torsion groups of exponential growth for which this property is established. Other examples that will be treated in this article include certain automorphisms groups of trees, $p$-adic general linear groups $GL_n(\mathbb{Q}_p)$ and special linear groups $SL_n(\mathbb{Q}_p)$, in their totally disconnected topology.

This article is organized as follows. In Section \ref{capacity subsection} we review Halmos' notion of capacity in Banach algebras and the Jenkins, Fountain, Ramsay and Williamson's result relating capacity and spectrum of selfadjoint elements. In Section \ref{growth subsection} we develop the appropriate notions and results regarding growth rates in locally compact groups.

Our main result and its corollaries, which relate fast growth with non-Hermitianess, are stated and proven in Section \ref{main result section}.

Finally, in Section \ref{applications section}, we give some known and also some new examples of non-Hermitian groups, using the methods developed in Section \ref{main result section}. We also state some open questions in Subsection \ref{questions subsection}.

\section{Preliminaries}
\label{preliminaries section}

Given a $^*$-algebra $A$ and an element $a \in A$ we will use throughout this article the notations $\sigma(a)$ and $R(a)$ to denote, respectively, the spectrum and the spectral radius of $a$. Also, if $S$ is a set, $|S|$ stands for the cardinality of $S$.\\

\subsection{Capacity in Banach algebras}
\label{capacity subsection}

In \cite{halmos} Halmos introduced the notion of \emph{capacity} of an element of a Banach algebra, so as to  give an appropriate analytic generalization, for Banach algebras, of the notion of an algebraic element (much like topological nilpotent elements are an analytic generalization of nilpotent elements). Halmos' definition, which can be found in \cite[page 857]{halmos}, was the following:\\

\begin{df}[Halmos]
Let $A$ be a Banach algebra. The \emph{capacity} of an element $a \in A$, denoted $\mathrm{cap}(a)$, is the number defined in the following way:
\begin{align}
\label{capacity formula def}
\mathrm{cap}(a) := \lim_{n} \inf_{p \in \mathcal{M}_{n}} \| p(a) \|^{\frac{1}{n}} = \lim_{n} \inf_{p \in \mathcal{P}_{n-1}} \| a^n + p(a) \|^{\frac{1}{n}}\,,
\end{align}
where $\mathcal{M}_n$ is the set of all monic complex polynomials of degree $n$ and $\mathcal{P}_n$ is the set of all complex polynomials of degree $n$. The limit in (\ref{capacity formula def}) can be shown to exist always.\\
\end{df}

Halmos' definition of capacity had roots in the classical notion of \emph{capacity} of a subset $X \subseteq \mathbb{C}$ coming from potential theory. The relation between the two concepts was established by Halmos himself when he showed that the capacity of an element $a$ of a Banach algebra is the same as the classical notion of capacity of the set $\sigma(a)$. This highlights the fact that the capacity of an element is in some way related with certain properties of its spectrum. In this regard, it is for instance clear  that the capacity of an element $a$ is always bounded by its spectral radius. Another important result in this setting is due Jenkins \cite[Corollary 1.4]{jenkins nonsym} and, apparently independently, to Fountain, Ramsay and Williamson \cite[Lema 5.1]{fountain}, which says that for a selfadjoint element in a Banach $^*$-algebra to have real spectrum, its capacity must necessarily be less than half of its spectral radius:\\

\begin{thm}[Jenkins \cite{jenkins nonsym}, Fountain, Ramsay, Williamson \cite{fountain}]
\label{fountain ramsay williamson theorem}
Let $A$ be a Banach $^*$-algebra and $a\in A$ a selfadjoint element. If $\sigma(a) \subseteq \mathbb{R}$, then $\mathrm{cap}(a) \leq \frac{1}{2}R(a)$.\\
\end{thm}

Jenkins did not state the above result in this way, but in an equivalent form. The above formulation appears in \cite{fountain}, and a proof of the result can be found there as well.

Theorem \ref{fountain ramsay williamson theorem} is a useful tool for showing that certain Banach $^*$-algebras are not Hermitian: one just needs to find a self-adjoint element for which $\mathrm{cap}(a) > \frac{1}{2}R(a)$. Fountain, Ramsay and Williamson used this result to give a new proof that the free group in two generators is not Hermitian \cite[pages 246-247]{fountain} and to give the first example of a non-Hermitian locally finite group \cite[page 248]{fountain}. The above result was also used later by Bomash \cite{bomash} to show that certain solvable groups are not Hermitian. Jenkins used the result to prove that $PGL_2(\mathbb{Q}_p)$ is not Hermitian in its totally disconnected topology \cite[page 300]{jenkins nonsym}.

In all these applications the authors established that a given element in the group algebra had a ``large'' capacity (greater than half of its spectral radius). In the case of Fountain, Ramsay and Williamson's new proof that the free group on two generators is not Hermitian, it is clear that the estimation of the capacity is based on the very fast growth of the free group. This is the idea we will follow in the remaining part of this article: we will show that if a subset of a group grows too fast, then there is a selfadjoint element with large capacity, and therefore the group cannot be Hermitian. We start, in the next subsection, by establishing the appropriate notions of growth in locally compact groups.\\

\subsection{Growth in locally compact groups}
\label{growth subsection}

We recall the definition of the \emph{growth function} and \emph{growth rate} of a subset of a locally compact group.\\

\begin{df}
Let $G$ be a locally compact group with a Haar measure $\mu$. For a measurable set $S \subseteq G$ the sequence $\{ \mu(S^n)^{\frac{1}{n}} \}_{n \in \mathbb{N}}$ is called the \emph{growth function} of $S$, and the limit superior 
\begin{align*}
\omega_{G}(S):= \limsup_{n \to \infty} \mu(S^n)^{\frac{1}{n}}
\end{align*}
is called the \emph{growth rate} of $S$.\\
\end{df}

Another important notion is that of \emph{spherical growth}:\\

\begin{df}
Let $G$ be a locally compact group. For a measurable set $S \subseteq G$ the sequence $\{ \mu(S^{n} \backslash S^{n-1})^{\frac{1}{n}} \}_{n \in \mathbb{N}}$ is called the \emph{spherical growth function} of $S$, and the limit superior 
\begin{align*}
\sigma_{G}(S):= \limsup_{n \to \infty} \mu(S^{n} \setminus S^{n-1})^{\frac{1}{n}}
\end{align*}
is called the \emph{spherical growth rate} of $S$.\\
\end{df}

It is clear that both the growth rate and the spherical growth rate do not depend on the choice of Haar measure $\mu$.

We have the following result:\\

\begin{thm}
\label{growth rate limit exists finite greater than one}
Let $G$ be a locally compact group and $S \subseteq G$ any nonempty relatively compact open set. We have that the limit
\begin{align*}
\lim_{n \to \infty} \mu(S^n)^{\frac{1}{n}}\,,
\end{align*}
always exists, is finite, and is greater or equal to $1$.\\
\end{thm}

The above result was established for compactly generated groups by Guivarc'h \cite[Th\' eoreme I.1]{guiv}, in the case where $S$ is a compact neighbourhood of $e$ that generates $G$. Since we are interested in general locally compact groups (not necessarily compactly generated) and since it will be especially important for us to consider sets that do not contain the identity $e$, we were lead to state and prove the result in the generality provided above. Before we prove this result, we recall the following standard definition:\\
\begin{df}
 A relatively compact open subset $S \subseteq G$ is said to have \emph{exponential growth} if $1  < \sigma_G(S)$ and \emph{subexponential growth} if $\sigma_G(S) = 1$.\\
\end{df}
We will now turn to the proof of the above result.\\

{\bf \emph{Proof of Theorem \ref{growth rate limit exists finite greater than one}:}} By [Guivarch Lemme I.1] we have that
\begin{align}
\label{guivarch lemma inequality 1}
\mu(S)\mu(S^{n+k}) \leq \mu(S^n S) \mu(S^{-1}S^k) = \mu(S^{n+1}) \mu(S^{-1}S^k)\,,
\end{align}
for any $n, k \in \mathbb{N}$. Also by [Guivarch Lemme I.1] we have that
\begin{align}
\label{guivarch lemma inequality 2}
\mu(S^{-1})\mu(S^{-1}S^k) \leq \mu(S^{-1}S^{-1}) \mu(SS^k) = \mu(S^{-2}) \mu(S^{k+1})\,.
\end{align}
 Taking the decomposition $S^{m + 1} = S^{m-1}S^2$ for any $m > 2$ and applying [Guivarch Lemme I.1] one last time we have that
\begin{align}
\label{guivarch lemma inequality 3}
\mu(S)\mu(S^{m+1}) \leq \mu(S^{m-1}S) \mu(S^{-1}S^2) = \mu(S^m) \mu(S^{-1}S^2)\,.
\end{align}

Using inequalities (\ref{guivarch lemma inequality 1}), (\ref{guivarch lemma inequality 2}) and (\ref{guivarch lemma inequality 3}) and the fact that $S$ is nonempty we obtain:
\begin{eqnarray*}
\mu(S^{n+k}) & \leq & \frac{1}{\mu(S)} \mu(S^{n+1}) \mu(S^{-1}S^k)\\
& = & \frac{\mu(S^{-2})}{\mu(S)\mu(S^{-1})} \mu(S^{n+1}) \mu(S^{k+1})\\
& = & \frac{\mu(S^{-2})\mu(S^{-1}S^2)^2}{\mu(S)^3\mu(S^{-1})}\; \mu(S^{n}) \mu(S^{k})\,.
\end{eqnarray*}
Hence, the sequence $\{\log \big( \mu(S^n)^{\frac{1}{n}} \big)\}_{n \in \mathbb{N}}$ satisfies the conditions of [Guivarch Lemme I.2] and we therefore conclude that the limit $\lim_{n \to \infty} \mu(S^n)^{\frac{1}{n}}$ exists and is finite. It is clear that this limit is always greater or equal to $1$ since $0 < \mu(S) \leq \mu(S^n)$ for all $n \in \mathbb{N}$. \qed\\

The spherical growth also always satisfies $1 \leq \sigma_G(S) < \infty$, provided one is only considering relatively compact symmetric sets $S$ that do not generate a compact subgroup, as we will now see. If $S$ generates a compact subgroup then it can happen, for example, that $\sigma_G(S) = 0$. We will also see later in this section, in Proposition \ref{prop limsup is a lim}, that under appropriate assumptions the limsup in the definition of $\sigma_G(S)$ is also a true limit, as we showed to be true for $\omega_G(S)$.\\

\begin{prop}
\label{spherical growth is finite and greater than one}
Let $G$ be a locally compact group and $S \subseteq G$ a relatively compact symmetric open set that does not generate a compact subgroup of $G$. Then the spherical growth $\sigma_G(S)$ is always finite and greater or equal to $1$.\\
\end{prop}

{\bf \emph{Proof:}} It is clear that the spherical growth is finite since $\sigma_{G}(S) \leq \omega_G(S)$ and $\omega_G(S)$ is finite by Theorem \ref{growth rate limit exists finite greater than one}.

 Suppose that $\sigma_G(S) < 1$. Then there is a real number $r$ such that eventually we have $\mu(S^n \backslash S^{n-1})^{\frac{1}{n}} < r < 1$, say for $n \geq n_0$. This means that, for $n \geq n_0$, we have $\mu(S^n \backslash S^{n-1}) < r^n$ and therefore the series $\sum_{n \geq 2} \mu(S^n \backslash S^{n-1})$ converges.
 
The subgroup $\langle S \rangle \subseteq G$ generated by $S$ is precisely the set
\begin{align}
\label{form of the subgroup gen by S}
\langle S \rangle = \bigcup_{n \geq 1} S^n\,,
\end{align}
because $S$ is symmetric (notice, for instance, that $e \in S^2$). Since we naturally have
\begin{align*}
\bigcup_{n \geq 1} S^n = S \cup \bigcup_{n \geq 2} S^n \backslash S^{n-1}\,,
\end{align*}
we can conclude that $\mu(\langle S \rangle )$ has finite measure because
\begin{align*}
\mu(\langle S \rangle )\;\; =\;\; \mu \Big(S \cup \bigcup_{n \geq 2} S^n \backslash S^{n-1} \Big)\;\; =\;\; \mu(S) + \sum_{n \geq 2} \mu(S^n \backslash S^{n-1}) \;\;<\;\; \infty\,.
\end{align*}
 From (\ref{form of the subgroup gen by S}) it is clear that $\langle S \rangle$ is an open subgroup, and therefore it is also closed. Since it has finite measure, it must be compact, which is a contradiction. Hence, we conclude that $1 \leq \sigma_G(S)$. \qed\\

It is a well-known fact that for non-compact finitely generated discrete groups the growth rate and the spherical growth rate coincide. We will now show that this holds for general locally compact groups too. This result will be very useful for us in the remaining sections of this article.\\

\begin{thm}
\label{growth rate and spherical growth rate coincide}
Let $G$ be a locally compact group and $S \subseteq G$ a symmetric relatively compact open set that does not generate a compact subgroup. We have that the growth rate and the spherical growth rate of $S$ coincide, i.e.
\begin{align*}
\omega_{G}(S) = \sigma_{G}(S)\,.
\end{align*}
\end{thm}

{\bf \emph{Proof:}} In this proof we consider two cases. First we consider the case where the sequence $\big(\mu(S^n \backslash S^{n-1})\big)_{n \in \mathbb{N}}$ is bounded by a constant $C \in \mathbb{R}^+$. In this case we necessarily have that
\begin{align*}
1 \;\leq\; \limsup_{n \to \infty} \mu(S^n \backslash S^{n-1})^{\frac{1}{n}} \;\leq\; \limsup_{n \to \infty} C^{\frac{1}{n}} \;= \; 1\,.
\end{align*}
So that $\sigma_G(S) = \limsup_{n \to \infty} \mu(S^n \backslash S^{n-1})^{\frac{1}{n}} = 1$. We also have that
\begin{eqnarray*}
 \omega_{G}(S) & = & \limsup_{n \to \infty} \mu(S^n )^{\frac{1}{n}} 
\;\; \leq \;\;  \limsup_{n \to \infty} \Big( \mu(S) +  \sum_{i = 1}^{n-1} \mu(S^{i+1} \backslash S^i) \Big)^{\frac{1}{n}}\\
& \leq & \limsup_{n \to \infty} \big( \mu(S)+ nC \big)^{\frac{1}{n}} \;\; = \;\;  1\,.
\end{eqnarray*}
Since it is always true that $1 \leq \omega_G(S)$, we conclude that $\omega_G(S) = 1$, and therefore $\omega_G(S) = \sigma_G(S)$.

We now consider the case where the sequence $\big(\mu(S^n \backslash S^{n-1})\big)_{n \in \mathbb{N}}$ is not bounded. It is clear that $\omega_G(S) \geq \sigma_G(S)$ so we only need to prove that $\omega_G(S) \leq \sigma_G(S)$.

 Since $\big(\mu(S^n \backslash S^{n-1})\big)_{n \in \mathbb{N}}$ is not bounded we can always find a subsequence $\big(\mu(S^{n_k} \backslash S^{n_{k}-1})\big)_{k \in \mathbb{N}}$ such that $ \mu(S^{i} \backslash S^{i-1}) \leq \mu(S^{n_k} \backslash S^{n_{k}-1})$ for any $i \leq n_k$. By Theorem \ref{growth rate limit exists finite greater than one} we have
\begin{eqnarray*}
\omega_G(S) & = & \lim_{n \to \infty} \mu(S^n )^{\frac{1}{n}} \;\; = \;\; \lim_{k \to \infty} \mu(S^{n_k} )^{\frac{1}{n_k}}\,.
\end{eqnarray*}
By the choice of the subsequence $\big(\mu(S^{n_k} \backslash S^{n_{k}-1})\big)_{k \in \mathbb{N}}$ we then have
\begin{eqnarray*}
\omega_G(S) & \leq &  \limsup_{k \to \infty} \Big( \mu(S) +  \sum_{i = 1}^{n_k-1} \mu(S^{i+1} \backslash S^i) \Big)^{\frac{1}{n_k}}\\
& \leq & \limsup_{k \to \infty}  \big(\mu(S) + n_k\, \mu(S^{n_k} \backslash S^{n_k -1}) \big)^{\frac{1}{n_k}}\\
& \leq & \limsup_{k \to \infty}  \big(n_k\, \mu(S^{n_k} \backslash S^{n_k -1}) \big)^{\frac{1}{n_k}}\\
& \leq & \big(\limsup_{k \to \infty}  n_k^{\frac{1}{n_k}}\big)\big(\limsup_{k \to \infty} \mu(S^{n_k} \backslash S^{n_k -1})^{\frac{1}{n_k}} \big)\\
& \leq & \limsup_{n \to \infty} \mu(S^{n} \backslash S^{n-1})^{\frac{1}{n}}\,.
\end{eqnarray*}
This finishes the proof. \qed\\

We now turn to the fact that, under certain assumptions on the set $S$, the limsup in the definition of the spherical growth rate $\sigma_G(S)$ is in fact a true limit:\\

\begin{prop}
\label{prop limsup is a lim}
Let $G$ be a locally compact group and $S$ a symmetric relatively compact open subset of $G$. Let us assume the following two conditions:
\begin{itemize}
\item[i)] $S$ does not generate a compact subgroup of $G$.
\item[ii)]  There exists a measurable symmetric set $K \subseteq G$ with $0 < \mu(K)$ and such that $SK \subseteq S$.
\end{itemize}
Then the limit $\lim_{n \to \infty} \mu(S^n \backslash S^{n-1})^{\frac{1}{n}}$ exists and is greater or equal to one.\\
\end{prop}

The above result is well-known in the case of discrete groups (where condition $ii)$ is trivially satisfied by taking $K= \{e\}$). We do not known if the result holds for a general subset $S$, i.e. if the result still holds without assuming condition $ii)$. Nevertheless, as we will see at the end of the section, this condition is satisfied by many subsets of totally disconnected locally compact groups.

In order to prove Proposition \ref{prop limsup is a lim} we will need the following lemma:\\

\begin{lemma}
\label{Sn - Sn-1 K = Sn - Sn-1}
Let $G$ be a locally compact group and $S$ a relatively compact open subset of $G$. If $K \subseteq G$ is a symmetric set such that $SK \subseteq S$, then for every $n \in \mathbb{N}$ we have $\big(S^n \backslash S^{n-1} \big)K \subseteq \big(S^n \backslash S^{n-1} \big)$.\\
\end{lemma}

{\bf \emph{Proof:}} If $S^n \backslash S^{n-1} = \emptyset$ the result is obvious. Let us assume then that $S^n \backslash S^{n-1} \neq \emptyset$ and let $g \in S^n \backslash S^{n-1}$ and $k \in K$. Since $SK \subseteq S$ we have $S^n K \subseteq S^n$, and therefore $gk \in S^n$. We now want to prove that $gk \notin S^{n-1}$. Suppose by contradiction that $gk \in S^{n-1}$. Then $g \in S^{n-1}k^{-1}$ and since $K$ is symmetric we have $g \in S^{n-1} K \subseteq S^{n-1}$, which is a contradiction. Thus, $\big(S^n \backslash S^{n-1} \big)K \subseteq \big(S^n \backslash S^{n-1} \big)$ \qed\\

{\bf \emph{Proof or Proposition \ref{prop limsup is a lim}:}} We claim that condition $i)$ implies that the sets $S^{n+1} \backslash S^{n}$ are non-empty for every $n \in \mathbb{N}$. If this was not the case, and $S^{n+1} \backslash S^{n} = \emptyset$ for a certain $n \geq 1$, then $S^{n+1} \subseteq S^n$ and consequently $S^k \subseteq S^n$ for every $k > n$. This would imply that the subgroup generated by $S$, which is $\bigcup_{n \in \mathbb{N}} S^n$, could be covered by the sets $S^1, \dots , S^n$, and would be therefore compact. Since this contradicts $i)$, we then know that the sets $S^{n+1} \backslash S^{n}$ are all non-empty.

 We claim that $S^{n+k} \backslash S^{n+k-1} \subseteq \big(S^{n} \backslash S^{n-1} \big) \big(S^{k} \backslash S^{k-1} \big)$ for every $n, k \in \mathbb{N}$. To see this, let $s_1, \dots ,s_n,t_1, \dots, t_k \in S$ be such that
\begin{align*}
s_1 \dots s_nt_1 \dots t_k \in S^{n+k} \backslash S^{n+k-1}\,.
\end{align*}
Then $s_1 \dots s_n \in S^n \backslash S^{n-1}$, because if $s_1 \dots s_n \in S^{n-1}$ then we would have $s_1 \dots s_nt_1 \dots t_k \in S^{n+k-1}$, which is not true by assumption. The same reasoning shows that $t_1 \dots t_k \in S^{k} \backslash S^{k-1}$. Hence $s_1 \dots s_nt_1 \dots t_k \in \big(S^{n} \backslash S^{n-1} \big) \big(S^{k} \backslash S^{k-1} \big)$.

By [Guivarch Lemme I.1] we have that
\begin{align*}
\mu(K) \mu\Big( \big(S^{n} \backslash S^{n-1} \big) \big(S^{k} \backslash S^{k-1} \big) \Big) \leq \mu\big(\big(S^{n} \backslash S^{n-1} \big) K \big) \mu\big( K^{-1}\big(S^{k} \backslash S^{k-1} \big)\big)\,.
\end{align*}

By Lemma \ref{Sn - Sn-1 K = Sn - Sn-1} we know that $\big(S^{n} \backslash S^{n-1} \big) K \subseteq \big( S^{n} \backslash S^{n-1} \big)$, and similarly, using the symmetry of $S$, we know that $K^{-1}\big(S^{n} \backslash S^{n-1} \big) \subseteq \big( S^{n} \backslash S^{n-1} \big)$. We can then conclude that
\begin{align*}
\mu\Big( \big(S^{n} \backslash S^{n-1} \big) \big(S^{k} \backslash S^{k-1} \big) \Big) \leq \frac{1}{\mu(K)} \mu\big(S^{n} \backslash S^{n-1} \big) \mu\big(S^{k} \backslash S^{k-1} \big)\,.
\end{align*}

As was shown before, we have that $S^{n+k} \backslash S^{n+k-1} \subseteq \big(S^{n} \backslash S^{n-1} \big) \big(S^{k} \backslash S^{k-1} \big)$, so that
\begin{align*}
\mu\big(S^{n+k} \backslash S^{n+k-1} \big) \leq \frac{1}{\mu(K)} \mu\big(S^{n} \backslash S^{n-1} \big) \mu\big(S^{k} \backslash S^{k-1} \big)\,.
\end{align*}

Hence, by [Guivarch Lemme I.2], the limit $\lim_{n \to \infty} \mu(S^n \backslash S^{n-1})^{\frac{1}{n}}$ must exist. The fact that this limit is finite and greater or equal to $1$ was already shown in Proposition \ref{spherical growth is finite and greater than one}. \qed\\

As we stated before, condition $ii)$ in Proposition \ref{prop limsup is a lim} is very suitable when dealing with totally disconnected groups, as the following result (which we are sure to be folklore for the experts in totally disconnected groups) shows.\\

\begin{prop}
\label{existence of K in tot disc groups}
 Let $G$ be a locally compact totally disconnected group. Let $S$ be any compact open subset of $G$. There exists a compact symmetric open subset $K \subseteq G$ such that $SK = S$.\\
\end{prop}

{\bf \emph{Proof:}} For every $g \in S$ choose a compact open subgroup $H_g$ inside $g^{-1}S$. We have that $S = \bigcup_{g \in S} gH_g$. Since $S$ is compact we can extract a finite cover $g_1H_{g_1}, \dots, g_n H_{g_n}$ of $S$. Take then any compact open subgroup $K$ such that $K \subseteq \bigcap_{i = 1}^n H_{g_i}$. Thus, since for any $i = 1, \dots, n$ we have that $K$ is a subgroup of $H_{g_i}$, we necessarily have $H_{g_i} K = H_{g_i}$. Hence, we conclude that
\begin{eqnarray*}
 SK & = & \bigcup_{i = 1}^n g_i H_{g_i} K \;\; = \;\; \bigcup_{i = 1}^n g_i H_{g_i} \;\; = \;\; S\,.
\end{eqnarray*}\qed\\

\section{Main Result}
\label{main result section}

The following theorem is the main result of this article and is our tool for establishing that certain groups are not Hermitian based on how fast some of their subsets grow.\\

\begin{thm}
\label{main theorem}
 Let $G$ be a locally compact group. Suppose there exists a relatively compact symmetric set $S \subseteq G$ that does not generate a compact subgroup and a relatively compact set $K \subseteq G$ such that $SK \subseteq S$. If we have
\begin{align}
\label{inequality for omegaS}
 \omega_G(S) > \frac{\int_S \Delta^{-\frac{1}{2}}(h) \;d\mu(h)}{2\, \mu(K)\,\inf_{h \in S} \Delta^{-\frac{1}{2}}(h)}\,,
\end{align}
 then $G$ is not Hermitian. Particularly, the element $f \in L^1(G)$ given by
\begin{align*}
 f(g):=\frac{1}{\int_S \Delta^{-\frac{1}{2}}(h) \;d\mu(h)}\;\; \Delta^{-\frac{1}{2}}(g)\chi_S(g)\,,
\end{align*}
 is a selfadjoint element with non-real spectrum.\\
\end{thm}

As we saw at the end of the previous section (Proposition \ref{existence of K in tot disc groups}), the condition on the existence of a relatively compact set $K \subseteq G$ such that $SK = S$ is naturally satisfied for totally disconnected locally compact groups. Hence, Theorem \ref{main theorem} is more easily applied in the totally disconnected setting and all the examples we will provide in Section \ref{applications section} belong to this class of groups.\\

\begin{lemma}
\label{inf Delta on Sn is geq than inf Delta exp n}
 Let $S$ be a subset of $G$ and $n \in \mathbb{N}$. We have that
\begin{align*}
 \inf_{g \in S^n} \Delta^{-\frac{1}{2}}(g) \geq \big( \inf_{g \in S} \Delta^{-\frac{1}{2}}(g) \big)^n\,.
\end{align*}\\
\end{lemma}

{\bf \emph{Proof:}} We have that
\begin{eqnarray*}
 \inf_{g \in S^n} \Delta^{-\frac{1}{2}}(g) & = &  \inf_{\substack{g_1 \in S \\ \vdots \\ g_n \in S}} \Delta^{-\frac{1}{2}}(g_1 \dots g_n) \;\; = \;\; \inf_{\substack{g_1 \in S \\ \vdots \\ g_n \in S}} \Delta^{-\frac{1}{2}}(g_1) \dots  \Delta^{-\frac{1}{2}}(g_n)\\
 & \geq & \big(\inf_{g_1 \in S} \Delta^{-\frac{1}{2}}(g_1) \big) \dots  \big(\inf_{g_n \in S}\Delta^{-\frac{1}{2}}(g_n) \big)\\
 & = & \big( \inf_{g \in S} \Delta^{-\frac{1}{2}}(g) \big)^n\,.
\end{eqnarray*}\qed\\

{\bf Proof of Theorem \ref{main theorem}:} To check that $f$ is selfadjoint we see that
\begin{eqnarray*}
 f^*(g) & = & \frac{1}{\int_S \Delta^{-\frac{1}{2}}(h) \;d\mu(h)}\;\; (\Delta^{-\frac{1}{2}}\chi_S)^*(g) \\
& = & \frac{1}{\int_S \Delta^{-\frac{1}{2}}(h) \;d\mu(h)}\;\; \Delta(g^{-1})\Delta^{-\frac{1}{2}}(g^{-1})\overline{\chi_S(g^{-1})}\\
& = & \frac{1}{\int_S \Delta^{-\frac{1}{2}}(h) \;d\mu(h)}\;\; \Delta^{-\frac{1}{2}}(g)\overline{\chi_S(g^{-1})}\\
& = & \frac{1}{\int_S \Delta^{-\frac{1}{2}}(h) \;d\mu(h)}\;\; \Delta^{-\frac{1}{2}}(g)\chi_S(g)\\
& = & f(g)\,.
\end{eqnarray*}

 The spectral radius of $f$ is less than $1$, because
\begin{eqnarray*}
 \rho(f) & \leq & \|f\|_1 \;\; = \;\;\frac{1}{\int_S \Delta^{-\frac{1}{2}}(h) \;d\mu(h) } \int_G \Delta^{-\frac{1}{2}}(g) \chi_S(g)\;d\mu(g)\\
& = & \frac{\int_S \Delta^{-\frac{1}{2}}(g) \;d\mu(g)}{\int_S \Delta^{-\frac{1}{2}}(h) \;d\mu(h)} \;\; = \;\;1\,.
\end{eqnarray*}
Hence by the Jenkins, Fountain, Ramsay, Williamson theorem (Theorem \ref{fountain ramsay williamson theorem} of the present article), we only need to prove that $\mathrm{cap}(f) > 1/ 2$ in order to show that $f$ has non-real spectrum, and hence that $G$ is not Hermitian.

We claim that for any $n \geq 1$ we have, for any $g \in S^n$,
\begin{align}
\label{inequality fn(g)}
 f^n (g) \geq \frac{\Delta^{-\frac{1}{2}}(g) \, \mu(K)^{n-1} }{\big( \int_S \Delta^{-\frac{1}{2}}(h) \;d\mu(h) \,\big)^n}\,.
\end{align}

 We will prove this claim by induction on $n$. The case $n = 1$ follows easily from the definition of $f$. Let us now consider the case $n \Rightarrow n +1$. Assume that the inequality (\ref{inequality fn(g)}) holds for a certain $n$. Any $g \in S^{n+1}$ can be written as $g = g_1 \dots g_{n+1}$, with $g_1, \dots, g_{n+1} \in S$. We have
\begin{eqnarray*}
 f^{n+1}(g) & = & \frac{1}{\int_S \Delta^{-\frac{1}{2}}(h) \;d\mu(h)} \int_G f^n (h) \Delta^{-\frac{1}{2}}(h^{-1}g)\chi_{S}(h^{-1}g)\;d\mu(h)\\
& \geq & \frac{1}{\int_S \Delta^{-\frac{1}{2}}(h) \;d\mu(h)} \int_{g_1\dots g_nK} f^n (h) \Delta^{-\frac{1}{2}}(h^{-1}g)\chi_{S}(h^{-1}g)\;d\mu(h)\,.
\end{eqnarray*}
Since $SK \subseteq S$, it follows that $S^n K \subseteq S^n$, and therefore $g_1\dots g_n K \subseteq S^n$. Hence, by the induction hypothesis, we have that
\begin{eqnarray*}
 f^{n+1}(g) & \geq & \frac{\mu(K)^{n-1}}{\big(\int_S \Delta^{-\frac{1}{2}}(h) \;d\mu(h) \,\big)^{n+1}} \int_{g_1\dots g_nK} \Delta^{-\frac{1}{2}} (h) \Delta^{-\frac{1}{2}}(h^{-1}g)\chi_{S}(h^{-1}g)\;d\mu(h)\\
& = & \frac{\Delta^{-\frac{1}{2}} (g)\,\mu(K)^{n-1}}{\big(\int_S \Delta^{-\frac{1}{2}}(h) \;d\mu(h) \,\big)^{n+1}} \int_{g_1\dots g_nK}  \chi_{S}(h^{-1}g)\;d\mu(h)\\
& = & \frac{\Delta^{-\frac{1}{2}} (g)\,\mu(K)^{n-1}}{\big(\int_S \Delta^{-\frac{1}{2}}(h) \;d\mu(h) \,\big)^{n+1}} \int_{K}  \chi_{S}(h^{-1}g_n^{-1} \dots g_1^{-1}g)\;d\mu(h)\\
& = & \frac{\Delta^{-\frac{1}{2}} (g)\,\mu(K)^{n-1}}{\big(\int_S \Delta^{-\frac{1}{2}}(h) \;d\mu(h) \,\big)^{n+1}} \int_{K}  \chi_{S}(h^{-1}g_{n+1})\;d\mu(h)\,.
\end{eqnarray*}
Since $K^{-1}S \subseteq S$ it follows that $h^{-1}g_{n+1} \in S$, for any $h \in K$. Thus we conclude that
\begin{eqnarray*}
 f^{n+1}(g) & \geq &  \frac{\Delta^{-\frac{1}{2}} (g)\,\mu(K)^{n-1}}{\big(\int_S \Delta^{-\frac{1}{2}}(h) \;d\mu(h) \,\big)^{n+1}} \int_{K}  1\;d\mu(h)\\
& = & \frac{\Delta^{-\frac{1}{2}} (g)\,\mu(K)^{n}}{\big(\int_S \Delta^{-\frac{1}{2}}(h) \;d\mu(h) \,\big)^{n+1}}\,,
\end{eqnarray*}
and this proves our claim.

We now claim that
\begin{align*}
 \| f^n |_{S^n \backslash S^{n-1}} \|_1 \geq \frac{\big( \inf_{h \in S} \Delta^{-\frac{1}{2}}(h)\, \big)^n \mu(K)^{n-1} \mu(S^n \backslash S^{n-1})}{\big(\int_S \Delta^{-\frac{1}{2}}(h) \;d\mu(h) \,\big)^n}\,.
\end{align*}
This inequality follows from the previous claim and Lemma \ref{inf Delta on Sn is geq than inf Delta exp n} according to the following computation:
\begin{eqnarray*}
 \| f^n |_{S^n \backslash S^{n-1}} \|_1 & = & \int_{S^n \backslash S^{n-1}} f^n(h) \;d\mu(h)\\
& \geq & \frac{ \mu(K)^{n-1} }{\big( \int_S \Delta^{-\frac{1}{2}}(h) \;d\mu(h) \,\big)^n} \int_{S^n \backslash S^{n-1}} \Delta^{-\frac{1}{2}}(h) \;d\mu(h)\\
& \geq & \frac{ \mu(K)^{n-1} }{\big( \int_S \Delta^{-\frac{1}{2}}(h) \;d\mu(h) \,\big)^n} \big(\inf_{h \in S^n \backslash S^{n-1}} \Delta^{-\frac{1}{2}}(h) \big) \mu(S^n \backslash S^{n-1})\\
& \geq & \frac{ \big(\inf_{h \in S^n} \Delta^{-\frac{1}{2}}(h) \big) \mu(K)^{n-1} \mu(S^n \backslash S^{n-1})}{\big( \int_S \Delta^{-\frac{1}{2}}(h) \;d\mu(h) \,\big)^n} \\
& \geq & \frac{ \big(\inf_{h \in S} \Delta^{-\frac{1}{2}}(h) \big)^n \mu(K)^{n-1} \mu(S^n \backslash S^{n-1})}{\big( \int_S \Delta^{-\frac{1}{2}}(h) \;d\mu(h) \,\big)^n}\,.
\end{eqnarray*}

 We will now estimate the capacity of $f$.
\begin{align*}
 \mathrm{cap}(f) \; = \; \lim_{n} \inf_{p \in \mathcal{P}_{n-1}} \| f^n + p(f)\|_1^{\frac{1}{n}} \; = \; \lim_{n} \inf_{p \in \mathcal{P}_{n-1}} \| f^n|_{S^n \backslash S^{n-1}} + f^n|_{S^{n-1}} + p(f)\|_1^{\frac{1}{n}}\,.
\end{align*}
Since $\mathrm{supp}(p(f)) \subseteq S^{n-1}$ for any $p \in \mathcal{P}_{n-1}$, it follows that
\begin{eqnarray*}
 \mathrm{cap}(f) & = & \lim_n \inf_{p \in \mathcal{P}_{n-1}} \Big(\| f^n|_{S^n \backslash S^{n-1}} \|_1 + \|f^n|_{S^{n-1}} + p(f)\|_1\Big)^{\frac{1}{n}}\\
& \geq & \lim_n \inf_{p \in \mathcal{P}_{n-1}} \| f^n|_{S^n \backslash S^{n-1}} \|_1^{\frac{1}{n}}\\
& = & \lim_n \| f^n|_{S^n \backslash S^{n-1}} \|_1^{\frac{1}{n}} \\
& \geq & \lim_n \left(\frac{ \big(\inf_{h \in S} \Delta^{-\frac{1}{2}}(h) \big)^n \mu(K)^{n-1} \mu(S^n \backslash S^{n-1})}{\big( \int_S \Delta^{-\frac{1}{2}}(h) \;d\mu(h) \,\big)^n}\right)^{\frac{1}{n}} \\
& = &  \lim_n \frac{ \big(\inf_{h \in S} \Delta^{-\frac{1}{2}}(h) \big) \mu(K)^{(n-1) / n} \mu(S^n \backslash S^{n-1})^{\frac{1}{n}}}{ \int_S \Delta^{-\frac{1}{2}}(h) \;d\mu(h)}\\
& = & \frac{ \big(\inf_{h \in S} \Delta^{-\frac{1}{2}}(h) \big) \mu(K) \omega_{G}(S)}{ \int_S \Delta^{-\frac{1}{2}}(h) \;d\mu(h)}\,,
\end{eqnarray*}
where in the last step we used Theorem \ref{growth rate and spherical growth rate coincide}. Hence, since we are assuming that inequality (\ref{inequality for omegaS}) holds, we conclude that $\mathrm{cap}(f) > \frac{1}{2}$, and therefore the element $f$ has non-real spectrum. \qed\\

For discrete groups the statement of Theorem \ref{main theorem} can be greatly simplified:\\

\begin{cor}
\label{corollary discrete groups}
 Suppose $G$ is discrete. If there exists a finite symmetric set $S \subseteq G$ with $|S| \geq 2$ and for which the inequality
\begin{align*}
 \omega_G(S) > \frac{|S|}{2}\,,
\end{align*}
holds, then $G$ is not Hermitian. In particular, the function $f \in \ell^1(G)$ defined by
\begin{align*}
 f := \frac{1}{|S|}\chi_S\,,
\end{align*}
is a selfadjoint element with non-real spectrum.\\
\end{cor}

{\bf \emph{Proof:}} If $|S| \geq 2$ and $\omega_G(S) > \frac{|S|}{2}$ it follows immediately that $\omega_G(S) > 1$ and therefore $S$ cannot generate a finite subgroup. The result then follows immediately from Theorem \ref{main theorem} by taking the counting measure as the Haar measure in $G$ and by taking $K$ as the trivial subgroup $\{e\}$.\qed\\

The following corollary of Theorem \ref{main theorem} is especially useful for tackling the question of Hermitianess in certain totally disconnected groups, as we shall see in Section \ref{applications section}.\\

\begin{cor}
\label{corollary totally disc selfadj groups}
 Let $G$ be a locally compact group and $K \subseteq G$ a compact open subgroup. We choose the normalized Haar measure $\mu$ such that $\mu(K) = 1$. If there exists a $g \in G$ satisfying the following conditions
\begin{itemize}
\item[i)] $KgK = Kg^{-1}K$,
\item[ii)] $\mu(KgK) > 1$,
\item[iii)] $\omega_G(KgK) > \frac{\mu(KgK)}{2}$,
\end{itemize} 
 then $G$ is not Hermitian.\\
\end{cor}

{\bf \emph{Proof:}} We take $S := KgK$, which is clearly a compact open symmetric set satisfying $SK = S$.  Since $\mu(KgK) \geq 2$ and $\omega_G(KgK) > \frac{\mu(KgK)}{2}$ it immediately follows that $S$ has exponential growth, and therefore does not generate a compact subgroup of $G$.

The modular function $\Delta$ is constant on the double coset $KgK$, since $K$ is a compact subgroup. Hence we have
\begin{align*}
\frac{\int_{KgK} \Delta^{-\frac{1}{2}}(h) \;d\mu(h)}{2\, \mu(K)\,\inf_{h \in KgK} \Delta^{-\frac{1}{2}}(h)} = \frac{\mu(KgK)}{2}\,.
\end{align*}
From $iii)$ we conclude that
\begin{align*}
\omega_G(S)  >  \frac{\int_{S} \Delta^{-\frac{1}{2}}(h) \;d\mu(h)}{2\, \mu(K)\,\inf_{h \in S} \Delta^{-\frac{1}{2}}(h)}\,,
\end{align*}
and it follows from Theorem \ref{main theorem} that $G$ is not Hermitian. \qed\\

\begin{rem}
\label{remark equality in main result}
A natural question to ask is if the growth condition (\ref{inequality for omegaS}) in Theorem \ref{main theorem} could be improved in order to ensure that $G$ is not Hermitian. The answer to this question is that, without any further assumptions, condition (\ref{inequality for omegaS}) is the sharpest possible, meaning that there are non-compact Hermitian groups $G$ for which
\begin{align*}
\omega_G(S) = \frac{\int_S \Delta^{-\frac{1}{2}}(h) \;d\mu(h)}{2\, \mu(K)\,\inf_{h \in S} \Delta^{-\frac{1}{2}}(h)}\,,
\end{align*}
for certain sets $S$ and $K$. A example of this is when $G = \mathbb{Z}$, $S =\{1, -1\}$ and $K = \{0\}$.\\
\end{rem}

\section{Applications}
\label{applications section}

In this section we will provide new examples of non-Hermitian groups, and also recover some known examples,  based on the results of the previous section that relate growth and Hermitianess.\\

\subsection{Some known results}
\label{some known results section}

\begin{cor}
Suppose $G$ is discrete. If $G$ contains a non-abelian free subgroup, then $G$ is not Hermitian.\\
\end{cor}

{\bf \emph{Proof:}} Let $a,b \in G$ be two elements that generate a non-abelian free subgroup of $G$. As it is well-known, the growth rate of $S := \{a, b, a^{-1}, b^{-1}\}$ is $\omega_G(S) = 3$. Thus we have $\omega_G(S) > \frac{|S|}{2} = 2$, and therefore $G$ is not Hermitian by Corollary \ref{corollary discrete groups}. \qed\\

It is known that, more generally, the existence of a free sub-semigroup in two generators is enough to ensure that $G$ is not Hermitian. It seems unlikely, however, that we can derive this more general result from Corollary \ref{corollary discrete groups} with the above proof, without any further assumptions. The problem lies with the fact that if $a$ and $b$ generate a free sub-semigroup, then we can only estimate the growth rate of $S:= \{a, b, a^{-1}, b^{-1}\}$ to be $\omega_G(S) \geq 2$. Hence, it is in principle possible that
\begin{align*}
\omega_G(S) = 2 = \frac{|S|}{2}\,,
\end{align*}
and this is not enough to ensure that $G$ is not Hermitian (see Remark \ref{remark equality in main result}).\\

\begin{ex}
The fundamental group $\Gamma_g$ of a closed orientable surface of genus $g \geq 2$ is not Hermitian. It is known that such groups contain free subgroups in two generators, and therefore cannot be Hermitian, but we will now present an alternative proof based on growth rates. The group $\Gamma_g$ has the presentation:
\begin{align*}
\Gamma_g := \langle a_1, b_1, \dots, a_g, b_g \; | \; a_1b_1a_1^{-1}b_1^{-1} \dots a_gb_ga_g^{-1}b_g^{-1} = 1 \rangle\,.
\end{align*}
Let us consider the symmetric set of generators
\begin{align*}
S_g := \{a_1, a_1^{-1}, b_1, b_1^{-1}, \dots ,  a_g, a_g^{-1}, b_g, b_g^{-1} \}\,.
\end{align*}
By \cite[VII Proposition 15]{harpe} it is known that $\omega_{\Gamma_g}(S) \geq 4g - 3$ for any set $S$ of generators of $\Gamma_g$ (this means that $\Gamma_g$ has what is usually called \emph{uniformly exponential growth}). Hence we see that
\begin{align*}
\omega_{\Gamma_g}(S) \geq 4g - 3 > \frac{4g}{2} = \frac{|S_g|}{2}\,,
\end{align*}
for all $g \geq 2$. Thus, by Corollary \ref{corollary discrete groups} it follows that $\Gamma_g$ is not Hermitian.\\
\end{ex}

\begin{ex}
The modular group $G:= (\mathbb{Z} / 2\mathbb{Z}) * (\mathbb{Z} / 3\mathbb{Z})$ is not Hermitian. This group contains a free subgroup in two generators (just like any non-trivial free product of groups with the exception of $(\mathbb{Z} / 2\mathbb{Z}) * (\mathbb{Z} / 2\mathbb{Z})$), and therefore cannot be Hermitian, but this can also be seen by looking at the growth rates.

Let us say that $a$ is the generator of $(\mathbb{Z} / 2\mathbb{Z})$ and $b$ is the generator of $(\mathbb{Z} / 3\mathbb{Z})$. We consider the symmetric set $S:= \{a, ab, b^2a \}$. As described in \cite[VI 7]{harpe}, we have $\omega_G(S) = \frac{\sqrt{5} + 1}{2}$ and therefore
\begin{align*}
\omega_G(S) = \frac{\sqrt{5} + 1}{2} > \frac{3}{2} = \frac{|S|}{2}\,,
\end{align*}
so that $G$ cannot be Hermitian by Corollary \ref{corollary discrete groups}.

It is interesting to note, however, that the argument does not work if one chooses the more natural symmetric set $S'=\{a, b, b^2\}$. In this case we have (by \cite[VI 7]{harpe}) that $\omega_G(S) = \sqrt{2}$ so that
\begin{align*}
\omega_G(S) = \sqrt{2} \not> \frac{3}{2} = \frac{|S'|}{2}\,.
\end{align*}
This difference of behaviour of the growth rate with respect to a chosen set of generators (both $S$ and $S'$ generate $G$) highlights the fact that one often has to make a careful choice for the set $S$ in order to able to apply Theorem \ref{main theorem} (see also question 5 in subsection \ref{questions subsection}).\\
\end{ex}

\begin{rem}
Concerning known results, it is also interesting to notice that Jenkins proof that $PGL_2(\mathbb{Q}_p)$ is not  Hermitian in its totally disconnected topology (\cite[page 300]{jenkins nonsym}) is done along similar lines as our Corollary \ref{corollary totally disc selfadj groups}, despite the fact that Jenkins does not explicitly refer to growth rates.\\
\end{rem}

\subsection{Free Burnside groups}

The property of having a free sub-semigroup in two generators was, essentially, the only general criterium in the literature for deciding that a discrete group was not Hermitian. Discrete groups without free sub-semigroups were therefore largely outside the scope of the existing theory, with the only examples of discrete groups without free sub-semigroups that were known to be non-Hermitian being the Fountain, Ramsay Williamson group \cite[page 248]{fountain} and Hulanicki's group \cite[Section 4]{hul4} (both being locally finite groups).

We will now give the first examples of non-Hermitian discrete groups of exponential growth which do not have free sub-semigroups. Recall that a free Burnside group $B(m,n)$, where $m, n \in \mathbb{N}$, is the (unique) group  generated by $m$ elements and such that every element $g \in B(m,n)$ satisfies the law $g^n = e$. These are obviously torsion groups (hence, without free sub-semigroups) and are sometimes infinite and of exponential growth (see \cite[Theorem 2.15]{adian}).\\

\begin{prop}
 The free Burnside groups $B(m, n)$ are not Hermitian for $m > 1$ and odd $n \geq 665$. \\
\end{prop}

{\bf \emph{Proof}} It is enough to prove this result for $m =2$, since a free Burnside group $B(m, n)$ with a number of generators $m$ larger than two always contains $B(2,n)$ and Hermitianess passes to open subgroups (\cite[Theorem 12.5.18 (a)]{palmer}).

The spherical growth function for $B(2,n)$ with odd $n \geq 665$ (with respect to the canonical set of generators, which we denote by $S$) was estimated by Adian in his proof that the free Burnside groups have exponential growth. In the proof of \cite[Theorem 2.15]{adian} he shows that this spherical growth function is very close to that of a free group in two generators, in the sense that
\begin{align*}
|S^k \setminus S^{k-1}| \geq 4 \cdot (2.9)^{k-1}\,.
\end{align*} 
Thus, if $a_1, a_2$ are the canonical generators of $B(2, n)$, the exponential growth rate of $S = \{a_1, a_2, a_1^{-1}, a_2^{-1} \}$ in $B(2, n)$ is
\begin{align*}
\omega_{B(2, n)}(S) = \sigma_{B(2,n)}(S) = \lim_{k} |S^k \setminus S^{k-1}|^{\frac{1}{n}} \geq 2.9\,.
\end{align*}
Hence, we clearly have $\omega_{B(2, n)}(S) > \frac{|S|}{2} = 2$, and by Corollary \ref{corollary discrete groups} we conclude that $B(2, n)$ is not Hermitian. \qed\\

\subsection{Automorphism groups of trees}

We start by recalling some of the terminology and facts about trees and their automorphisms that we are going to use. Recall that a \emph{tree} is a connected graph with no cycles. A tree $X$ has a natural \emph{distance function} $d:X^0 \times X^0 \to \mathbb{N}$ on the set of its vertices $X^0$.  By a \emph{geodesic} $L$ in a tree we mean a subgraph which is isomorphic to the real line $\mathbb{R}$. Also, the number of edges which are incident to a vertex $x_0$ is called the \emph{degree} of $x_0$, and the tree is said to be \emph{locally finite} if every vertex has finite degree.

 A \emph{tree automorphism} is a bijective map of the set of vertices of the tree onto itself which preserves the edges, and a tree automorphism is called a \emph{translation} of length $k \in \mathbb{N}$ if it has an (unique) invariant geodesic whose vertices are translated by a distance of length $k$.
 
 The group $Aut(X)$ of automorphisms of a locally finite tree $X$ can be naturally given a locally compact totally disconnected topology, with the sets $U_F(g) := \{h \in Aut(X) : g(x) = h(x)\,, \forall x \in F \}$, where $F$ is a finite set of vertices, forming a basis of neighbourhoods of the element $g \in Aut(X)$. In this topology, the stabilizer of a vertex is a compact open subgroup.
 
 We have the following result:\\

\begin{prop}
\label{automorphisms tree result}
Let $X$ be a locally finite tree and $G \subseteq Aut(X)$ a closed subgroup of automorphisms of $X$. Suppose that the following conditions are satisfied:
\begin{itemize}
\item[i)] $G$ contains a translation $g$ along a geodesic $L$, with translation length $k$;
\item[ii)] $L$ contains a point $x_0$ with degree at least $3$;
\item[iii)] The stabilizer $K$ of $x_0$ inside $G$ acts transitively on each sphere $\mathscr{S}_{n}(x_0) :=\{x \in X: d(x_0, x) = n\}$, with $n \in \mathbb{N}$.
\end{itemize}
Then $G$ is not Hermitian in its totally disconnected topology.\\
\end{prop}

{\bf \emph{Proof:}} Let $\mu$ be the normalized Haar measure of $G$ for which $\mu(K) =1$. Since the action of $K$ on the sphere $\mathscr{S}_{k}(x_0)$ is transitive and $g$ takes $x_0$ into a point in this sphere, it follows from \cite[Lemma 2.1 (3)]{hecke trees} that $KgK = Kg^{-1}K$.

It follows also from \cite[Lemma 2.1 (3)]{hecke trees} that $\mu(KgK)$ equals the cardinality of the sphere $\mathscr{S}_k(x_0)$. Since $x_0$ has degree greater or equal to $3$,  $X$ has a geodesic that contains $x_0$ and the action of $K$ is transitive on each sphere, there must be at least $3$ elements in $\mathscr{S}_k(x_0)$, and therefore $\mu(KgK) >1$.

Since the action of $K$ on each sphere $\mathscr{S}_{n}(x_0)$, where $n \in \mathbb{N}$, is transitive,  all the vertices of $\mathscr{S}_{n}(x_0)$ must have the same degree, which we denote by $\deg_n \in \mathbb{N}$. It is now not difficult to see that, for $n \geq 1$,
\begin{align*}
|\mathscr{S}_{n+1}(x_0)| = (\deg_n -1) \cdot  |\mathscr{S}_{n}(x_0)|\,,
\end{align*}
and of course $|\mathscr{S}_{1}(x_0)| = \deg_0$. Hence, the cardinality of the each sphere is given by 
\begin{align*}
|\mathscr{S}_{n}(x_0)| = \deg_0 \prod_{i = 1}^{n-1} (\deg_i - 1) \,.
\end{align*}
Since $g$ is a translation of translation length $k$, it is not difficult to see that $\deg_{n + k} = \deg_{n}$, for all $n \in \mathbb{N}$. Hence we have that
\begin{align*}
|\mathscr{S}_{nk}(x_0)| = \deg_0 (\deg_0 - 1)^{n-1}\prod_{i = 1}^{k-1} (\deg_i - 1)^n\,.
\end{align*}

Using \cite[Lemma 2.1 (3)]{hecke trees} again, we see that $\mu(Kg^nK)$ is the cardinality of the sphere $\mathscr{S}_{nk}(x_0)$. Hence, we see that
\begin{eqnarray*}
\mu((KgK)^n)  &\geq & \mu(Kg^nK) \;\; = \;\; |\mathscr{S}_{nk}(x_0)|\\
& = &  \deg_0 (\deg_0 - 1)^{n-1}\prod_{i = 1}^{k-1} (\deg_i - 1)^n \\
& \geq &   (\deg_0 - 1)^n\prod_{i = 1}^{k-1} (\deg_i - 1)^n\,.
\end{eqnarray*}
Observing that for any integer $q \geq 3$ we have $(q-1) \geq \frac{2}{3} q$, it follows that
\begin{eqnarray*}
\mu((KgK)^n) & \geq & \Big(\frac{2}{3}\deg_0 \Big)^n \prod_{i = 1}^{k-1} (\deg_i - 1)^n\\
& = & \big(\frac{2}{3} \big)^n\mu(KgK)^n\,.
\end{eqnarray*}
Hence we have 
\begin{align*}
\lim_{n \to \infty} \mu((KgK)^n)^{\frac{1}{n}}  \;\; =  \;\; \frac{2}{3} \mu(KgK) \;\; > \;\; \frac{1}{2}\mu(KgK)\,.
\end{align*}
The conditions of Corollary \ref{corollary totally disc selfadj groups} are then satisfied, and therefore $G$ cannot be Hermitian. \qed\\

Among the most studied trees in the literature are the so-called \emph{semihomogeneous} trees. These are the locally finite trees for which we can divide the set of vertices into two disjoint sets $X_1$ and $X_2$, with the vertices  of $X_1$ (respectively, $X_2$) all having the same degree, and such that every vertex of $X_1$ is connected only to vertices of $X_2$ and vice-versa. The automorphism group of such trees always has translations and the stabilizer of any vertex acts transitively on every sphere around it. Hence, we immediately have the following result:\\

\begin{cor}
Let $X$ be a semihomogeneous tree that has a vertex with degree at least $3$. Then its automorphism group $Aut(X)$ is not a Hermitian group in its totally disconnected topology.\\
\end{cor}

\begin{rem}
It also follows from Proposition \ref{automorphisms tree result} that the group $SL_2(\mathbb{Q}_p)$, where $\mathbb{Q}_p$ is the field of $p$-adic numbers, is not Hermitian in its totally disconnected topology, by considering it as a group of automorphisms of its associated Bruhat-Tits tree. This group was shown to be non-Hermitian in \cite[page 19]{palma}, using a different method.\\
\end{rem}

\subsection{General linear groups $GL_n(\mathbb{Q}_p)$}
 
Let $p$ be a prime number and $\mathbb{Q}_p$ be the field of $p$-adic numbers. 
 
Given an $n$-tuple $\lambda:=(\lambda_1, \dots, \lambda_n) \in \mathbb{Z}^n$ we will denote by $\pi^{\lambda}$ the diagonal matrix
\begin{align*} \pi^{\lambda} = \left(
\begin{array}{ccc}
p^{\lambda_1} & & \\
 & \ddots & \\
 & & p^{\lambda_n}
\end{array} \right) \in GL_n(\mathbb{Q}_p)\,.
\end{align*}
As it is known $K:= GL_n(\mathbb{Z}_p)$ is a compact open subgroup of $GL_n(\mathbb{Q}_p)$. Let $\mu$ denote the Haar measure of $GL_n(\mathbb{Q}_p)$, normalized so that $\mu(K) = 1$. In our next result we will need know to the precise value of measure of $K \pi^{\lambda} K$. The case where $\lambda \in \mathbb{N}^{n}$ has been computed in \cite[chapter V (2.9)]{sym} and we have that
\begin{align}
\label{expression for left cosets in GLn}
\mu(K\pi^{\lambda}K) = p^{2\langle \lambda \,,\, \rho \rangle} \nu_n(p^{-1})/ \nu_{\lambda}(p^{-1})\,,
\end{align}
where $\langle \cdot \,,\, \cdot \rangle$ is the usual inner product, $\rho := \frac{1}{2}(n-1, n-3, n-5, \dots, 1-n) \in \mathbb{Z}^n$, each function $\nu_m(t)$, defined in \cite[chapter III - 1]{sym},  is given by
\begin{align*}
\nu_m(t) = (1-t)^{-m}\prod_{i=1}^m 1-t^i\,,
\end{align*}
and for any $\lambda = (\lambda_1, \dots, \lambda_k)$ the function $\nu_{\lambda}(t)$, defined in \cite[chapter III - 1]{sym}, is given by
\begin{align*}
\nu_{\lambda}(t) = \prod_{i=1}^m \nu_{m_i}(t)\,,
\end{align*}
where $m_i$ is the number of $\lambda_j$ equal to $i$, for each $i \geq 0$.

We have the following result:\\

\begin{prop}
\label{GL is not Herm}
For any $n \geq 2$ and any prime $p \geq 5$, the group $GL_n(\mathbb{Q}_p)$ is not Hermitian in its totally disconnected topology.\\
\end{prop}

{\bf \emph{Proof:}}  Let $\lambda := (1, 0, 0, \dots, 0,0,-1)$ and $K:= GL_n(\mathbb{Q}_p)$. Since $K$ contains all the permutation matrices, it is clear that $K(\pi^{\lambda})^{-1}K = K \pi^{-\lambda} K = K\pi^{\lambda}K$.

For any $k \in \mathbb{N}$ it is clear that $(\pi^{\lambda})^k = \pi^{k\lambda}$. Our goal is to use Corollary \ref{corollary totally disc selfadj groups} and for that we have to compute the measure of $K\pi^{k\lambda}K$. Expression (\ref{expression for left cosets in GLn}) is only valid for $\lambda \in \mathbb{N}^n$, which is not the case here, but there is nevertheless a simple trick that allows us to reduce to this case: letting $\overline{k} \in \mathbb{N}^n$ be the constant tuple equal to $k$, we observe that
\begin{eqnarray*}
\mu(K\pi^{k\lambda}K) & = & \mu(\pi^{\overline{k}} K\pi^{k\lambda}K) \;\; = \;\; \mu( K\pi^{\overline{k}}\pi^{k\lambda}K)\\
& = & \mu( K\pi^{k\lambda + \overline{k}}K)\\
& = & p^{2\langle k\lambda + \overline{k} \,,\, \rho \rangle} \nu_n(p^{-1})/ \nu_{k\lambda + \overline{k}}(p^{-1})\,.
\end{eqnarray*}
It is not difficult to see that $\langle \overline{k} \,,\, \rho \rangle = 0$, so that
\begin{align*}
p^{2\langle k \lambda + \overline{k}\,,\, \rho \rangle} = p^{2\langle k \lambda\,,\, \rho \rangle} = p^{k(n-1) - k(1-n)} = p^{2k(n-1)}\,.
\end{align*}
Moreover, since there is precisely one entry in $k\lambda + \overline{k}$ that equals $0$, one entry that equals $2k$, and $n-2$ entries that equal $k$, an easy computation yields that $\nu_{k\lambda + \overline{k}}(p^{-1}) = \nu_1(p^{-1})^2\nu_{n-2}(p^{-1})$. Hence, we have that
\begin{eqnarray*}
\mu(K\pi^{k\lambda}K)& = & p^{2k(n-1)} \frac{\nu_n(p^{-1})}{\nu_1(p^{-1})^2\nu_{n-2}(p^{-1})}\\
& = & p^{2k(n-1)} \frac{(1-p^{-(n-1)})(1-p^{-n})}{(1-p^{-1})^2}\,.
\end{eqnarray*}
From this, it readily follows that $\mu(K \pi^{\lambda} K) > 1$ and moreover
\begin{align*}
\lim_{k \to \infty} \mu((K\pi^{\lambda}K)^k)^{\frac{1}{k}} \geq \lim_{k \to \infty} \mu(K\pi^{k\lambda}K)^{\frac{1}{k}} = p^{2(n-1)}\,.
\end{align*}
In order for us to use Corollary \ref{corollary totally disc selfadj groups} we will then need to show that  $p^{2(n-1)} > \frac{\mu(K\pi^{\lambda}K)}{2}$, which amounts to showing that
\begin{align*}
p^{2(n-1)} > \frac{p^{2(n-1)}(1-p^{-(n-1)})(1-p^{-n})}{2(1-p^{-1})^2}\,,
\end{align*}
or equivalently
\begin{align}
\label{inequality in proof of GLn}
2(1-p^{-1})^2 - (1-p^{-(n-1)})(1-p^{-n}) > 0 \,.
\end{align}
We will prove that this is true for any $n \geq 2$ and any prime $p \geq 5$. To see this we consider the function $f(t) = 2(1-t^{-1})^2 - 1$. We clearly have that
\begin{align*}
2(1-t^{-1})^2 - (1-t^{-(n-1)})(1-t^{-n}) > f(t) \,,
\end{align*}
for $t \geq 0$, so we just need to see that $f(t) > 0$ for all $t \geq 5$. We have that $f(5) = 2(\frac{4}{5})^2 - 1 =  \frac{32}{25} - 1 > 0$. Moreover, we have that $f'(t) = 4(1-t^{-1})t^{-2}$, so that $f'(t) > 0$ for $t \geq 5$. Hence, $f$ is surely growing from the point $t = 5$ onwards, and therefore $f(t) >0$ for all $t \geq 5$.

By Corollary  \ref{corollary totally disc selfadj groups} we conclude that $GL_n(\mathbb{Q}_p)$ is not Hermitian for $p \geq 5$ and $n \geq 2$. \qed\\

\begin{rem}
It can also be seen from the proof above that $GL_2(\mathbb{Q}_2)$, $GL_2(\mathbb{Q}_3)$ and $GL_3(\mathbb{Q}_3)$ are not Hermitian in their totally disconnected topology, simply by checking that inequality (\ref{inequality in proof of GLn}) is satisfied with respect to these choices of $n$ and $p$.\\
\end{rem}

\begin{rem}
There are other choices of matrices $\pi^{\lambda}$ which could be usefully considered. For example, if $n$ is even, we could take $\lambda := (1, \dots, 1, 0, \dots, 0)$, where half of the entries equal $1$ and the other half equal $0$. If one works with $PGL(\mathbb{Q}_p)$ instead, we have $K [\pi^{\lambda}] K = K [\pi^{  \lambda}]^{-1} K$, where $K$ is the image of $GL_n(\mathbb{Z}_p)$ in $PGL_n(\mathbb{Q}_p)$. It would be possible to then use similar methods as above to prove, for the prime number $3$, that $PGL_n(\mathbb{Q}_3)$ is not Hermitian (hence, $GL_n(\mathbb{Q}_3)$ is not Hermitian) for various choices of $n$ as an even number.

We think therefore that it is reasonable to conjecture that the groups $GL_n(\mathbb{Q}_p)$ are not Hermitian for all $n \geq 2$ and all primes $p$.
\end{rem}

\subsection{Special linear groups $SL_n(\mathbb{Q}_p)$}

The argument we presented in the previous subsection to prove that $GL_n(\mathbb{Q}_p)$ is not Hermitian works as well for the special linear groups $SL_n(\mathbb{Q}_p)$.\\

\begin{prop}
\label{SL is not Herm}
For any $n \geq 2$ and any prime $p \geq 5$, the group $SL_n(\mathbb{Q}_p)$ is not Hermitian in its totally disconnected topology.\\
\end{prop}

{\bf \emph{Proof:}} Let $\mu$ and $\widetilde{\mu}$ be the normalized Haar measures of $SL_n(\mathbb{Q}_p)$ and $GL_n(\mathbb{Q}_p)$, respectively, for which the compact open subgroups $K := SL_n(\mathbb{Z}_p)$ and  $\widetilde{K} := GL_n(\mathbb{Z}_p)$ have measure $1$. We claim that if $g \in SL_n(\mathbb{Q}_p) \subseteq GL_n(\mathbb{Q}_p)$, then
\begin{align}
\label{claim about measures SL and GL}
\mu(KgK) = \widetilde{\mu}(\widetilde{K} g \widetilde{K})\,.
\end{align}
To prove this, we start by noticing that, since $K$ is a compact open subgroup, $\mu(KgK)$ is equal to the total number of left cosets inside $KgK$. Similarly, $\widetilde{\mu}(\widetilde{K} g \widetilde{K})$ is equal to the total number of left cosets inside $\widetilde{K} g \widetilde{K}$, so that we only need to prove that these numbers are the same.

Let $KgK / K$ and $\widetilde{K} g \widetilde{K} / \widetilde{K}$ be the sets of left cosets inside $KgK$ and $\widetilde{K} g \widetilde{K}$, respectively. We consider the following map:
\begin{align}
KgK / K & \longrightarrow \widetilde{K} g \widetilde{K} / \widetilde{K} \notag\\
hK & \mapsto h\widetilde{K}\,.
\label{map proof SL}
\end{align}
This map is easily seen to be well-defined. We claim that it is surjective. Let $kg \widetilde{K} \in \widetilde{K} g \widetilde{K} / \widetilde{K}$, with $k \in \widetilde{K}$. We have $\big(\det(k)^{-1}k \big)gK \in KgK / K$, and the image of this element via the map (\ref{map proof SL}) is precisely $kg\widetilde{K}$ because
\begin{align*}
\big(\det(k)^{-1}k \big)g\widetilde{K} \;=\; \big(\det(k)^{-1}k \big)g\big(\det(k)k^{-1} \big) \widetilde{K} \;=\; kgk^{-1} \widetilde{K} \;=\; kg\widetilde{K}\,.
\end{align*}
This proves that the map (\ref{map proof SL}) is surjective. Let us now prove that it is injective.
Let $k_1gK$ and $k_2gK$ be two left cosets inside $KgK$, with $k_1, k_2 \in K$. Suppose $k_1g\widetilde{K} = k_2g\widetilde{K}$. Then, we have $g^{-1}k_2^{-1}k_1g \in \widetilde{K}$. Since $\det(g^{-1}k_2^{-1}k_1g) = 1$ we must have $g^{-1}k_2^{-1}k_1g \in K$, and therefore $k_1gK = k_2gK$. Thus, injectivity of the map (\ref{map proof SL}) is proven, and this yields equality (\ref{claim about measures SL and GL}).

Given that the measures $\mu(KgK)$ and $\widetilde{\mu}(\widetilde{K} g \widetilde{K})$ are the same, the argument in the proof of Proposition \ref{GL is not Herm} used for showing that $GL_n(\mathbb{Q})$ is not Hermitian can also be applied to the group $SL_n(\mathbb{Q}_p)$, because the matrix $\pi^{\lambda}$ we start with (as well as its powers) is in $SL_n(\mathbb{Q}_p)$. \qed\\

\subsection{Some questions}
\label{questions subsection}

\begin{itemize}
\item[1)] Is the $p$-adic "$ax + b$" group $\mathbb{Q}_p \rtimes \mathbb{Q}_p^*$ Hermitian as a totally disconnected group?

It is known that $\mathbb{R} \rtimes \mathbb{R}^*$ is Hermitian as a connected Lie group \cite{leptin}, and $\mathbb{Q} \rtimes \mathbb{Q}^*$ is not Hermitian as a discrete group (it has free sub-semigroups). When trying to use the methods developed in this article to tackle the problem for  $\mathbb{Q}_p \rtimes \mathbb{Q}_p^*$, we have only been able to obtain expression (\ref{inequality for omegaS}) as an equality, which is not sufficient to assure the group is not Hermitian.

\item[2)] Suppose $G$ is a totally disconnected locally compact group, $H \subseteq G$ is a compact open subgroup and $\alpha$ an automorphism of $G$ such that $\alpha(H) \subsetneq H$ and $G = \bigcup_{n \in \mathbb{Z}} \alpha^n(H)$. Is the semi-direct product $G \rtimes_{\alpha} \mathbb{Z}$ Hermitian?

An easy example of such a group is $\mathbb{Q}_p \rtimes \mathbb{Z}$, where the automorphism is given by multiplication by $p$. Similarly to the group in question $1)$, we could only obtain expression (\ref{inequality for omegaS}) as an equality, which is insufficient to prove that $G \rtimes_{\alpha} \mathbb{Z}$ is not Hermitian.

\item[3)]  Are there totally disconnected Hermitian groups of exponential growth? Can such examples be found among discrete groups?

These questions have already been asked by the author in \cite{palma}, and the groups in questions $1)$ and $2)$ are natural to be considered in this regard. A negative answer to this question would imply that all Hermitian totally disconnected groups are amenable, which would give more evidence to the following long standing conjecture (see \cite{palmer}):

\item[4)] Are all Hermitian groups amenable?

This conjecture was answered affirmatively for connected groups by Palmer (see \cite[Theorem 12.5.18 (e)]{palmer}), but remains open in general, even for discrete groups.

\item[5)] Suppose $G$ is a finitely generated group, with $S$ being a symmetric finite set of generators. We  know that the growth rate of $S$ lies in between $1$ and the growth rate of a free group in $|S|$ generators, i.e. $1 \leq \omega_G(S) \leq |S| -1$.

Let us now normalize this value so that it becomes independent of the number of generators, i.e. let us consider the number $\theta_G(S)$ defined by
\begin{align*}
\theta_G(S) := \frac{\omega_G(S) - 1}{|S| -2} \in [0,1]\,.
\end{align*}
What can we say about $\sup_{S} \theta_G(S)$, where the supremum runs over all finite symmetric sets $S$ of generators?

If $\theta_G(S) = 0$, then $G$ has subexponential growth, while if $\theta_G(S) = 1$ then $G$ is necessarily a free group. The value of $\inf_S \theta_G(S)$ (without our normalization) has been widely studied, and is behind what is known as \emph{uniform exponential growth}. Understanding the supremum $\sup_{S} \theta_G(S)$ would give valuable information regarding the Hermitianess of $G$, as our Corollary \ref{corollary discrete groups} shows. 

\end{itemize}

\end{document}